\magnification=1200
\centerline{\bf A bound for the number of automorphisms of an arithmetic Riemann
surface}
\bigskip
\centerline{BY MIKHAIL BELOLIPETSKY}
\medskip
\centerline{{\it Sobolev Institute of Mathematics, Novosibirsk,} 630090,
{\it Russia}}
\centerline{\it e-mail: \tt mbel@math.nsc.ru}
\bigskip
\centerline{AND GARETH A.~JONES}
\medskip
\centerline{\it Faculty of Mathematical Studies, University of Southampton,}
\centerline{{\it Southampton, SO\/}17 1{\it BJ}}
\centerline{\it e-mail: \tt gaj@maths.soton.ac.uk}

\bigskip\bigskip

\centerline{\it Abstract}

\medskip

We show that for every $g\geq 2$ there is a compact arithmetic Riemann
surface of genus $g$ with at least
$4(g-1)$ automorphisms, and that this lower bound is attained by infinitely
many genera, the smallest being $24$.

\medskip

\moveright163pt \vbox{\hrule width2.4cm height1pt}

\bigskip

\centerline{1. \it Introduction}

\medskip

Schwarz [{\bf 17}] proved that the automorphism group of a compact
Riemann surface of genus $g\geq 2$ is finite, and Hurwitz [{\bf 10}]
showed that its order is at most $84(g-1)$. This bound is sharp,
by which we mean that it is attained for infinitely many $g$, and the least
genus of such an extremal
surface is $3$. However, it is also well known that there are infinitely
many genera for which the
bound $84(g-1)$ is not attained. It therefore makes sense to consider the
maximal order $N(g)$ of the
group of automorphisms of any Riemann surface of genus $g$. Accola [{\bf
1}] and Maclachlan
[{\bf 14}] independently proved that $N(g)\geq 8(g+1)$. This bound is also
sharp, and according
to p.~93 of [{\bf 2}], Paul Hewitt has shown that the least genus attaining
it is $23$.
Thus we have the following sharp bounds for $N(g)$ with $g\geq 2$:
$$8(g+1) \leq N(g) \leq 84(g-1).$$

We now consider these bounds from an arithmetic point of view, defining
arithmetic Riemann surfaces to be those which are uniformized by arithmetic
Fuchsian groups. The motivation for this approach can be found in the
works of Borel, Margulis and various others on arithmetic groups.
Concerning Riemann surfaces with large groups of automorphisms, the
surprising fact, which can easily be seen, is that all the extremal
surfaces for Hurwitz's upper bound are arithmetic, whereas all the
extremal surfaces for the Accola-Maclachlan lower bound are non-arithmetic.
This raises the natural question: ``What can be said about the other two
bounds?"

The non-arithmetic analog of Hurwitz's upper bound, obtained by the first
author in [{\bf 3}], is $156(g-1)/7$; this bound is sharp, and the least genus
attaining it is $50$. The aim of the current paper is to obtain an arithmetic
analog of the Accola-Maclachlan lower bound, namely that for each $g\geq 2$
there is an arithmetic surface of genus $g$ with $4(g-1)$ automorphisms,
and that this bound is attained for infinitely many $g$, starting with $24$.

We now collect these results together: defining $N_{\rm ar}(g)$ and $N_{\rm
na}(g)$ to be the maximal orders of the automorphisms groups of the arithmetic
and non-arithmetic surfaces of genus $g$ respectively, for all sufficiently
large $g$ we have a system of sharp bounds
$$4(g-1) \leq N_{\rm ar}(g) \leq 84(g-1),$$
$$8(g+1) \leq N_{\rm na}(g) \leq {156\over 7}(g-1).$$

In Section 2 we recall the basic facts about Riemann surfaces and
arithmetic groups. Section 3 contains the proof of the $4(g-1)$ lower
bound, with a number of additional remarks. Finally, in Section 4 we
use our proof of the $4(g-1)$ bound to describe an infinite set of genera
for which the bound is attained, and to prove that the least genus
attaining the bound is $24$.

The authors would like to thank A.~D.~Mednykh for suggesting this problem,
and C.~Maclachlan for some very helpful discussions.

\bigskip

\centerline{2. \it Basic facts}

\medskip

In this section we recall some definitions and basic properties of Riemann
surfaces and arithmetic groups. For more information about Riemann surfaces
and Fuchsian groups see [{\bf 7}, {\bf 11}]. The basic references for
quaternion algebras and arithmetic groups are [{\bf 12}, {\bf 20}].

\medskip

{\it Definition} 2.1. A {\it Riemann surface\/} is a connected
one-dimensional complex analytic manifold. An {\it automorphism\/}
of a Riemann surface is an analytic mapping of the surface onto
itself.

\medskip

In this paper we shall consider only compact Riemann surfaces of genus $g\geq
2$. By the uniformization theorem [{\bf 7}, Ch.~IV] each such surface ${\cal
S}$ can be represented as the quotient space ${\cal H}/\Gamma_{\cal S}$, where
$\cal H$ is the hyperbolic plane and $\Gamma_{\cal S}$ is a cocompact
torsion-free discrete subgroup of the group ${\rm Isom}^+({\cal H}) =
PSL(2,{\bf R})$ of orientation-preserving isometries of $\cal H$. This group
$\Gamma_{\cal S}$, called the {\it surface group\/} corresponding to ${\cal
S}$, is unique up to conjugacy in $PSL(2,{\bf R})$ and is finitely generated.

Discrete subgroups of $PSL(2,{\bf R})$ are called {\it Fuchsian groups}. Each
cocompact Fuchsian group $\Gamma$ has a {\it signature\/}
$\sigma=(g;m_1,\ldots, m_k)$, where $g$ is a non-negative integer, equal to the
genus of ${\cal H}/\Gamma$, and each $m_j$ is an integer greater than $1$,
indicating a cone-point of order $m_j$ in ${\cal H}/\Gamma$. This signature
corresponds to the canonical presentation for $\Gamma$:
$$
\Gamma (g;m_1,\dots ,m_k)
=\langle\,\alpha_1,\beta_1,\dots,\alpha_g,\beta_g,\gamma_1,\dots,\gamma_k
\mid\prod_{i=1}^g[\alpha_i,\beta_i]\prod_{j=1}^k\gamma_j=1,\,\gamma_j^{m_j}=
1\,\rangle.
$$
If $g=0$ we shall omit $g$ from $\sigma$, and write $(m_1,\dots ,m_k)$.

We define $\mu(\Gamma)$ to be the hyperbolic measure of ${\cal H}/\Gamma$, it
can be expressed in terms of the signature:
$$
\mu(\Gamma)=\mu(g;m_1,\dots ,m_k)=2\pi\left(2g-2+\sum_{j=1}^k
\left(1-{1\over m_j}\right)\right).\eqno(1)
$$

By the Riemann-Hurwitz formula if $\Gamma'$ is a subgroup of index $n$ in
$\Gamma$ we have:
$$\mu(\Gamma')=n\cdot\mu(\Gamma).$$

The automorphisms of a Riemann surface ${\cal S}$ lift to
the isometries of $\cal H$ normalizing the surface group $\Gamma_{\cal S}$,
so ${\cal S}$ has automorphism group
$${\rm Aut}\,({\cal S}) \cong N(\Gamma_{\cal S})/\Gamma_{\cal S}$$
where $N(\Gamma_{\cal S})$ is the normalizer of $\Gamma_{\cal S}$ in
$PSL(2,{\bf R})$.

In our investigations we often need to construct a Riemann surface ${\cal S}$
with a given Fuchsian group $\Gamma$ normalizing its surface group
$\Gamma_{\cal S}$. In order to do this one has to find a torsion-free normal
subgroup of finite index in $\Gamma$, or equivalently to find an epimorphism
from $\Gamma$ onto some finite group $G$ with a torsion-free kernel. We call
such an epimorphism a {\it surface-kernel epimorphism\/}, or {\it SKE\/} for
short. In these circumstances ${\rm Aut}\,({\cal S})$ has a subgroup isomorphic
to $G$. It is known that in any Fuchsian group all elements of finite order are
conjugate to powers of the elliptic generators in a canonical presentation of
the group. Hence, in order to verify that a given epimorphism is a SKE, one has
only to check that the orders of these generators are preserved.

Now we introduce a special class of Riemann surfaces, which we call
arithmetic surfaces.

\medskip

{\it Definition\/} 2.2. (See [{\bf 4}, {\bf 12}, {\bf 15}, {\bf 18}, {\bf 20}].)
Let $A=({{a,b}\over{k}})$ be a quaternion algebra over a totally real number
field $k$, such that there is an isomorphism $\rho$ from $({{a,b}\over{\bf R}})$
to the matrix algebra $M_2({\bf R})$ and such that
$({{\sigma(a),\sigma(b)}\over{\bf R}}) \cong {\bf H}$ (Hamilton's quaternions)
for every non-identity Galois monomorphism $\sigma: k\to {\bf R}$. Let $\cal O$
be an order in $A$, and let ${\cal O}^1$ be the group of
elements of norm $1$ in ${\cal O}$. Then any subgroup of $PSL(2,{\bf R})$ which
is commensurable with the image in $PSL(2,{\bf R})$ of some such
$\rho({\cal O}^1)$ is called an {\it arithmetic Fuchsian group}.

\medskip

Arithmeticity is invariant under conjugation in $PSL(2,{\bf R})$, so the
following definition is valid:

\medskip

{\it Definition\/} 2.3. A Riemann surface is {\it arithmetic\/} if it is
uniformized by an
arithmetic Fuchsian group. All other Riemann surfaces are {\it nonarithmetic}.

\medskip

We finish this section with some examples of arithmetic Fuchsian groups
and Riemann surfaces.

\medskip

{\it Example} 2.1. {\it Triangle groups\/} are Fuchsian groups
which have signatures of the form $(m_1, m_2 ,m_3)$. Triangle groups with a
given
signature are conjugate in $PSL(2,{\bf R})$ (this fails for most other
signatures), so either all of them or none of them are arithmetic. Takeuchi
first proved
that there are only finitely many signatures of arithmetic triangle groups,
and gave the complete list of them in [{\bf 18}]; particularly important
examples
for us are the signatures $(2,3,7)$ and $(2,4,5)$.
In order to obtain this result Takeuchi
used an arithmeticity test which he introduced in the same paper.

\medskip

{\it Example} 2.2.
The orientation-preserving subgroup of the group generated by reflections
in the sides of a right-angled hyperbolic pentagon $\Pi$ is a Fuchsian
group $\Gamma$ of signature $(2,2,2,2,2)$.
If $\Pi$ can be subdivided
into $n$ congruent triangles, so that whenever two triangles have a common
side they are
symmetric with respect to that side, then $\Gamma$ is a subgroup of index
$n$ in the corresponding
triangle group. In particular, if this triangle group is arithmetic then so
is $\Gamma$. For
instance, one can barycentrically subdivide an equilateral right-angled
pentagon into $10$ triangles
with angles $\pi/2$, $\pi/4$ and $\pi/5$, and so obtain a
$(2,2,2,2,2)$-subgroup of
index $10$ in the arithmetic $(2,4,5)$-group. However, it is worth noting
that among the arithmetic groups of a given signature there may also be maximal
Fuchsian groups, and these can not be obtained as subgroups of arithmetic
triangle
groups.

\medskip

{\it Example} 2.3. All surfaces of genus $g$ with $84(g-1)$
automorphisms (such as Klein's quartic) are arithmetic, since they are
uniformized by finite index subgroups of the $(2,3,7)$ triangle group,
which is arithmetic.

\bigskip

\centerline{3. \it The main results}

\medskip

LEMMA 3.1. {\it Let $\{{\cal S}_g\}_{g\in{\cal G}}$ be an infinite
sequence of arithmetic surfaces of different genera $g$, such that for
each $g\in\cal G$ the group of automorphisms of ${\cal S}_g$ has order
$a(g+b)$ for some fixed $a$ and $b$. Then $b=-1$.}

\medskip

{\it Proof.} Let ${\cal S}$ be a surface from the given
sequence. Since ${\rm Aut}\,({\cal S})\cong N(\Gamma_{\cal S})/\Gamma_{\cal
S}$, the Riemann-Hurwitz formula gives
$$\mu(N(\Gamma_{\cal S}))={\mu(\Gamma_{\cal S})\over |{\rm Aut}\,({\cal S})|}
={2\pi(2g-2)\over a(g+b)},$$
so $\mu(N(\Gamma_{\cal S}))\to 4\pi/a$ as $g\to\infty$.

Since $\Gamma_{\cal S}$ is an arithmetic Fuchsian group, $N(\Gamma_{\cal
S})$ is also arithmetic. Borel [{\bf 4}] showed that the measures of
arithmetic groups form a discrete subset of $\bf R$, so for all but finitely
many $g\in{\cal G}$ we have
$${2\pi(2g-2)\over a(g+b)} = {4\pi\over a},$$
and from this it follows that $b=-1$.

\medskip

As an immediate consequence of Lemma~3.1 we deduce that the Accola-Maclachlan
lower bound for $N(g)$ cannot be attained by infinitely many arithmetic
surfaces. In fact, since the extremal surfaces for this bound are uniformized by
surface subgroups of $(2,4,2(g+1))$-groups with $g\geq 24$ [{\bf 14}], and
these are not arithmetic by [{\bf 18}], it is never attained by arithmetic
surfaces.

It also follows from Lemma~3.1 that the infinite sequences of Riemann
surfaces with automorphism groups of order $8(g+1)$, $8(g+3)$, etc., as
studied by Accola [{\bf 1}], Conder and Kulkarni [{\bf 5}], Maclachlan
[{\bf 14}], and others, can be constructed only in non-arithmetic
situations.

We now come to the central question of this paper, which is to find
a sharp lower bound for $N_{\rm ar}(g)$.

\medskip

LEMMA 3.2. {\it $N_{\rm ar}(g)\geq 4(g-1)$ for all $g\geq 2$.}

\medskip

{\it Proof.} Let
$\Gamma=\langle \gamma_1,\ldots,\gamma_5\mid \gamma_j^2=\gamma_1\ldots
\gamma_5=1\rangle$
be an arithmetic group with signature $(2,2,2,2,2)$ (see Example 2.2).
Consider the homomorphism $\theta$ from $\Gamma$ to the dihedral group
$G=D_{2(g-1)}=\langle a,b\mid a^{2(g-1)}=b^2=(ab)^2=1\rangle$ of
order $4(g-1)$ defined by $\gamma_j\mapsto ab,\, b,\, a^{g-2}b,\, b,\, a^{g-1}$
for $j=1,\ldots, 5$. It is easy to verify that $\theta$ is a SKE. The kernel
$K={\rm ker}\,(\theta)$ is therefore a surface group, and the surface
${\cal S} = {\cal H}/K$ is arithmetic since $K$ is a finite index subgroup
of the arithmetic group $\Gamma$. Since $\mu(\Gamma)=\pi$ and $|G|=4(g-1)$,
the Riemann-Hurwitz formula gives $\mu(K) = \mu(\Gamma)|G| = 2\pi(2g-2)$
and so ${\cal S}$ has genus $g$. Since ${\rm Aut}\,({\cal S})\geq \Gamma/K
\cong G$
it follows that $N_{\rm ar}(g)\geq |G|=4(g-1)$.

\medskip

THEOREM 3.1. {\it  $N_{\rm ar}(g)\geq 4(g-1)$ for all $g\geq 2$, and this
bound is
attained for infinitely many values of $g$.}

\medskip

{\it Proof.} The inequality in the statement of the theorem was proved in
the previous
lemma, so it remains to show that the bound is sharp. Suppose that $G:={\rm
Aut}\,({\cal S})$ has order $|G|>4(g-1)$ for some compact arithmetic
surface $\cal S$ of genus $g\geq
2$. We will successively impose a set of conditions on $g$ which lead to a
contradiction, and then
show that infinitely many values of $g$ satisfy these conditions.

By our hypothesis, $G\cong\Gamma/K$ for some cocompact
arithmetic group $\Gamma$ and normal surface subgroup $K=\Gamma_{\cal S}$
of $\Gamma$, with
$$4\pi(g-1)=\mu(K)=|G|\mu(\Gamma)>4(g-1)\mu(\Gamma), \eqno(2)$$
so $\mu(\Gamma)<\pi$. Borel's discreteness theorem [{\bf 4}] implies that
there are only
finitely many measures of cocompact arithmetic groups $\mu(\Gamma)<\pi$,
and then formula (1) for
$\mu(\Gamma)$ shows that these correspond to a finite set $\Sigma$ of
signatures, all of genus
$0$ and with either three or four elliptic periods.

For each $\sigma\in\Sigma$, the number $q=\mu(\Gamma)/4\pi$ is rational
and depends only on the signature $\sigma$ of $\Gamma$, so write
$q=r/s=r_{\sigma}/s_{\sigma}$
in reduced form, and let $R={\rm lcm}\{\,r_{\sigma}\mid
\sigma\in\Sigma\,\}$. By~$(2)$ we have $|G|=(g-1)/q=(g-1)s/r$ for some such
$q$. Since $|G|$ is an integer, if we choose $g$ so that $g-1$ is coprime
to $R$ then for surfaces of genus $g$ we have $r=1$ and $|G|=(g-1)s$.

Suppose that $g-1$ is a prime $p>S$, where
$S=\max\{\,s_{\sigma}\mid\sigma\in\Sigma,\,
r_{\sigma}=1\,\}$. Then $|G|=ps$ with $s$ coprime to
$p$ and less than $p+1$, so Sylow's Theorems imply that $G$ has a normal
Sylow $p$-subgroup $P\cong C_p$. Let $\Delta$ denote the inverse image of
$P$ in $\Gamma$, a normal subgroup of $\Gamma$ with $\Gamma/\Delta\cong
Q:=G/P$ of order $s$. Let $\Pi$ denote the finite set of primes which divide an
elliptic period $m_j$ of some signature
$\sigma\in\Sigma$ with $r_{\sigma}=1$. If we take $p\not\in\Pi$, then the
natural
epimorphism $G\to Q$ preserves the orders of the images of all elliptic
generators of $\Gamma$; the inclusion $K\leq\Delta$ therefore induces a
smooth $p$-sheeted covering ${\cal S}\to{\cal T}={\cal H}/\Delta$ of
surfaces, so $\Delta$ is a surface group of genus $1+(g-1)/p=2$. Thus $Q$
is a group of automorphisms of a Riemann surface ${\cal T}$ of genus $2$,
so $Q$ is one of a known finite list of groups (for instance, $|Q|\leq
48$). Let $E$ denote the least common multiple of the exponents of all
the groups of automorphisms of Riemann surfaces of genus $2$. (All we
need here is the fact that $E$ is finite and even.)

Since $\Delta/K\cong P\cong C_p$ it follows that $K$ contains the subgroup
$\Delta'\Delta^p$ generated by the commutators and $p$-th powers in $\Delta$,
so $P$ is isomorphic (as a ${\bf Z}_pQ$-module) to a 1-dimensional quotient
of the ${\bf Z}_pQ$-module $\Delta/\Delta'\Delta^p$, where the action of $Q$
is induced by conjugation in $\Gamma$. Now $\Delta$ is isomorphic to the
fundamental group $\pi_1({\cal T})$ of $\cal T$, so $\Delta/\Delta'$ is
isomorphic (as a ${\bf Z}Q$-module) to its first integer homology group
$H_1({\cal T},{\bf Z})\cong\pi_1({\cal T})^{\rm ab}$, and hence
$\Delta/\Delta'\Delta^p$ is isomorphic (as a ${\bf Z}_pQ$-module) to
$H_1({\cal T},{\bf Z}_p) \cong H_1({\cal T},{\bf Z})\otimes{\bf Z}_p$;
since $\cal T$ has genus $\gamma=2$, this has dimension $2\gamma=4$
as a vector space over ${\bf Z}_p$. Since $p$ does not divide
$s=|Q|$, Maschke's Theorem [{\bf 8}, I.17.7] implies that
$H_1({\cal T},{\bf Z}_p)$ is a direct sum of irreducible submodules.
Now $H_1({\cal T},{\bf C})=H_1({\cal T},{\bf Z})\otimes{\bf C}$ is a
direct sum of two $Q$-invariant subspaces, corresponding under duality
to the holomorphic and antiholomorphic differentials in $H^1({\cal T},
{\bf C})$, and these afford complex conjugate representations of $Q$
[{\bf 16}]. It follows that there are just three possibilities for
$H_1({\cal T},{\bf Z}_p)$: it may be irreducible, a direct sum of two
irreducible 2-dimensional submodules, or a direct sum of four irreducible
1-dimensional submodules. Since $H_1({\cal T},{\bf Z}_p)$ has a
$1$-dimensional quotient, only the last of these three cases can arise.
We have $p>2$ (since $p>S\geq 2$), so a theorem of Serre [{\bf 7}, V.3.4]
implies that $Q$ is faithfully represented on $H_1({\cal T},{\bf Z}_p)$;
thus $Q\leq GL_1(p)^4\cong(C_{p-1})^4$, so $Q$ has exponent $e$ dividing
$p-1$. Since $e$ also divides $E$, if we choose $p$ so that
$\gcd(p-1,E)=2$ then $e$ must divide $2$. Since $\Delta$ is a surface
group, the natural epimorphism $\Gamma\to\Gamma/\Delta\cong Q$ is a
SKE, so each elliptic period of $\Gamma$ is equal to $2$. However, as
noted earlier, $\Gamma$ is a cocompact Fuchsian group of genus $0$
with at most four elliptic periods, so this contradicts the fact that
$\mu(\Gamma)>0$.

It remains to check that there are infinitely many values of $g$
satisfying the above conditions, namely that $g-1$ is a prime $p$
where $p>S$, $p\not\in\Pi$, $p$ is coprime to $R$, and $\gcd(p-1,E)=2$.
Dirichlet's Theorem implies that there are infinitely many primes $p$
satisfying the last condition (for instance, primes $p\equiv -1$
mod~$(E)$), and all but finitely many satisfy the other three conditions
(since $\Pi$ is finite), so the proof is complete.

\medskip

{\it Remark} 3.1. In this proof arithmeticity is used only to
show that there are just finitely many signatures $\sigma$ that can occur.
It follows that there are similar results for other classes of groups
with this finiteness property.

\medskip

{\it Remark} 3.2. For our chosen values of $g$, the bound $4(g-1)$
is attained only by dihedral quotients of $\Gamma=\Gamma(2,2,2,2,2)$, as in
Lemma~3.2. To see this, repeat the proof of Theorem~3.1, but starting
with $|{\rm Aut}\,({\cal S})|\geq 4(g-1)$ instead of strict inequality. We
eventually find that $\Gamma=\Gamma(2,2,2,2,2)$ or $\Gamma=\Gamma(1;2)
=\langle \alpha,\beta,\gamma\mid \gamma^2=[\alpha,\beta]\gamma=1\rangle$;
since $Q$ is abelian (having exponent $2$), all commutators in $G$ lie in
$P$ and hence there is no SKE from $\Gamma(1;2)$ onto $Q$. Thus $\Gamma=
\Gamma(2,2,2,2,2)$. The Riemann-Hurwitz formula gives $|Q|=\mu(\Delta)/
\mu(\Gamma)=4\pi/\pi=4$, so $Q\cong V_4$ (a Klein four-group). Since
${\rm Aut}\,(C_p)\cong C_{p-1}$ the only extensions of $C_p$ by $V_4$ are
$C_p\times V_4$ and $D_p\times C_2\cong D_{2p}$; there is no epimorphism
$\Gamma(2,2,2,2,2)\to C_p$, so we must have $G\cong D_{2p}=D_{2(g-1)}$.

\medskip

{\it Remark} 3.3. Theorem~3.1 has another more elementary
proof. Let us start with the sequence of genera in [{\bf 14}] which
attain the Accola-Maclachlan bound. These have the form $g=89p+1$,
with the prime $p$ satisfying five additional conditions. As mentioned
earlier, the extremal surfaces of these genera are uniformized by surface
subgroups of the
$(2,4,2(g+1))$ triangle groups, which are non-arithmetic, so for such $g$
we have $4(g-1)\leq N_{\rm ar}(g)<N(g)=8(g+1)$. There are only finitely many
signatures $\sigma$ with $\mu(2,4,2(g+1))<\mu(\sigma)<\mu(2,2,2,2,2)$
which can correspond to arithmetic groups, and so may be considered as
candidates for giving a better lower bound for $N_{\rm ar}(g)$. Using the
Riemann-Hurwitz formula one can apply divisibility arguments to exclude
those signatures which do not have surface subgroups of genus $89p+1$.
Together with known information about arithmetic groups of signature
$(2,2,2,n)$ [{\bf 15}, {\bf 19}] and arithmetic triangle
groups [{\bf 18}], this gives the following list of candidates:
$(2,5,20)$, $(2,6,12)$, $(2,8,8)$, $(3,4,6)$, $(4,4,4)$, $(2,2,2,4)$,
$(2,7,14)$, $(2,9,18)$, $(2,12,12)$, $(3,4,12)$, $(3,6,6)$, $(4,4,6)$,
$(2,2,2,6)$, $(2,2,3,3)$, $(2,15,30)$, $(5,5,5)$, $(2,2,2,10)$.
Using a case-by-case argument, one can show if $p$ is a sufficiently large
prime then no group
$\Gamma$ of such a signature can have a normal surface subgroup of genus
$g=89p+1$. If we also
impose the conditions on $p$ given in [{\bf 14}] then we obtain a sequence
of genera $g$ for which the
arithmetic bound $4(g-1)$ is sharp.

In this approach one needs only Sylow's Theorems and some other basic
facts about finite groups, but the proof is rather routine and not very
straightforward: it is easy to handle the signatures with large elliptic
periods, but it becomes more complicated when the periods are small.
The most challenging case is when $\sigma=(2,2,2,4)$. The other reason
why we prefer our initial proof of Theorem~3.1 will be clear after
the next section, where we find the minimal genus for which our bound is
attained.

\bigskip

\centerline{4. \it Extremal surfaces}

\medskip

In this section we shall first use the proof of Theorem~3.1 to produce a
specific set of genera $g$ attaining our lower bound for $N_{\rm ar}(g)$. We
shall then strengthen the arguments in order to consider smaller $g$, and
finally determine the least genus for which $N_{\rm ar}(g)=4(g-1)$.

To begin with, let us see which signatures actually form the set $\Sigma$
corresponding to the cocompact arithmetic groups $\Gamma$ with
$\mu(\Gamma)<\pi$. Firstly, almost all of the cocompact arithmetic
triangle groups in Takeuchi's list [{\bf 18}] are contained in $\Sigma$.
Simple calculations show that the other possible signatures are $(2,2,2,n)$
for $n\geq 3$, $(2,2,3,3)$, $(2,2,3,4)$ and $(2,2,3,5)$. The arithmetic
groups of signature $(2,2,2,n)$ with odd $n$ were determined by Maclachlan
and Rosenberger [{\bf 15}]. For even $n$, groups of signature $(2,2,2,n)$
have a subgroup of index $2$ isomorphic to a $(1;n/2)$-group, and the
list of arithmetic groups of signature $(1;n/2)$ was obtained by Takeuchi
[{\bf 19}]. Combining these results we find that only $12$ signatures
of the form $(2,2,2,n)$ yield arithmetic groups. It is a matter of
direct verification whether or not there are arithmetic groups of the remaining
three signatures, but since this does not affect our arguments we shall
ignore this point and include them in the table of signatures $\sigma\in\Sigma$
given at the end of the paper.

Now inspecting this list $\Sigma$ of possible signatures, one can use the
proof of Theorem~3.1 to produce specific values of $g$ attaining the
lower bound $N_{\rm ar}(g)\geq 4(g-1)$. For instance, we see that $R=2^2 \cdot
3\cdot 5\cdot 7$, so a prime $p$ is coprime to $R$ provided $p>7$.
Inspection also shows that $\Pi=\{2,3,5,7\}$, so $p\not\in\Pi$
if and only if $p>7$. We also have
$S=84$, so we need $p>84$ for the proof of the Theorem to work
(though it can be adapted to apply to certain smaller primes, as we shall
see). A standard result [{\bf 7}, V.1.11] states that a Riemann surface
of genus $\gamma\geq 2$ has no automorphisms of prime order greater than
$2\gamma+1$, so taking $\gamma=2$ we see that $E$ is divisible only by
the primes $2,3$ and $5$; hence the condition $\gcd(p-1,E)=2$ is satisfied
by all odd $p$ such that $p-1$ is not divisible by $3,4$ or $5$, that is,
$p\equiv 23, 47$ or $59$ mod~$(60)$. It follows that for all such primes
$p>84$, our bound is attained by $g=p+1$. The smallest prime in this sequence
is $p=107$, giving $g=108$.

If we inspect $\Sigma$ more closely, and use other group-theoretic
techniques in addition to Sylow's Theorems, we can find smaller values
of $g$ attaining our bound. The basic idea is that, in order to show that
$G$ has a normal Sylow $p$-subgroup of order $p$, we replace the rather
crude sufficient condition $p>S$ with a more careful analysis of the
possibilities for a group of order $ps$. We use the fact that (according to
the Table)
the largest possible values of $s$ (for $r=1$) are $s=84,\, 48,\, 40,\,
36,\, 30$,
corresponding to the arithmetic groups of signatures $(2,3,7),\, (2,3,8),\,
(2,4,5),\, (2,3,9),\, (2,3,10)$ respectively, followed by $s=24$ corresponding
to $(2,3,12)$ and $(2,4,6)$, and then $s=21$ corresponding to $(2,3,14)$.

\medskip

{\it Example} 4.1. Let $p=59$, so $g=60$. We follow the proof of
Theorem~3.1, amending it where necessary for this particular prime $p$.
Since $p>7$, $p$ is coprime to $R$ and hence $|G|=(g-1)s=59s$. No possible
value of $s$ is divisible by $59$ (see the Table), so a Sylow $59$-subgroup
$P$ of $G$ has order $59$. The number of Sylow $59$-subgroups divides $s$
and is congruent to $1$ mod~$(59)$; this immediately implies that there is
only one such subgroup, so $P$ is normal in $G$. For the rest of the
proof, it is sufficient to note that $p=59$ satisfies $p\not\in\Pi$ and
$\gcd(p-1,E)=2$, so $g=60$ attains the lower bound.

\medskip

{\it Example} 4.2. Let $p=47$. Once again $p>7$, $p$ is coprime
to $R$, and $|G|=(g-1)s=47s$. A Sylow $47$-subgroup $P$ of $G$ has order
$47$ since there is no value of $s$ divisible by $47$. We need to show
that $P$ is normal in $G$, so suppose not. The number $n_{47}$ of Sylow
$47$-subgroups divides $s$ and is congruent to $1$ mod~$(47)$, so (by
inspection) $n_{47}=s=48$. This means that $P=N_G(P)$, so $G$ permutes
its $48$ Sylow $47$-subgroups by conjugation as a Frobenius group. A
theorem of Frobenius [{\bf 8}, V.7.6, V.8.2] implies that $G$ has a normal
subgroup $N$ of order $48$ (the Frobenius kernel), so $\Gamma$ has an
epimorphism onto $G/N\cong C_{47}$. However, $s=48$ implies that $\Gamma$
is the triangle group $\Gamma(2,3,8)$, so no such epimorphism exists, and hence
$P$ is normal in $G$. The rest of the proof is the same, so the lower bound
is attained for $g=48$. This method also deals with $g=84$, using
$p=83$ and $\Gamma=\Gamma(2,3,7)$.

\medskip

{\it Example} 4.3. Let $p=23$. As with $p=47$, the only place where
the proof of Theorem~3.1 fails is that Sylow's Theorems are not strong
enough to prove that a Sylow $23$-subgroup $P$ of $G$ is normal and has
order $23$. By inspection of $\Sigma$, no possible value of $s$ is
divisible by $23$, so $|P|=23$. Similarly, if $P$ is not normal, then
there must be $n_{23}=24$ Sylow $23$-subgroups, with $s=24$ or $48$, so
$|G|=24\cdot 23$ or $48\cdot 23$. In either case, $G$ permutes its Sylow
$23$-subgroups by conjugation as a transitive permutation group $\widetilde
G$ of degree
$24$. In fact, $\widetilde G$ is doubly transitive, since $P$ must have a
single orbit of length $23$ on the remaining Sylow $23$-subgroups: if it
normalized a Sylow $p$-subgroup $P^*\neq P$, then $PP^*$ would be a subgroup
of $G$ of order $23^2$. Thus $|\widetilde G|$ is divisible by $24\cdot 23$,
and it divides $|G|$, so $|\widetilde G|=24\cdot 23$ or $48\cdot 23$. In the
first case, $\widetilde G$ is sharply $2$-transitive, which is impossible
since such groups all have prime-power degree [{\bf 9}, XII.9.1]. In
the second case, since a point-stabilizer must act as $D_{23}$,
$\widetilde G$ is a Zassenhaus group with two-point stabilizers of even
order ($=2$); Zassenhaus showed that such groups of degree $n$ have
two-point stabilizers of order at least $(n-2)/2$ [{\bf 9}, XI.1.10],
so this case is also impossible. (Alternatively, the classification of
finite simple groups implies that the doubly transitive finite groups are all
known [{\bf 6}]: those of degree $24$ are the symmetric group $S_{24}$, the
alternating group $A_{24}$, the Mathieu group $M_{24}$, the projective
general linear group $PGL(2,23)$, and the projective special linear group
$PSL(2,23)$, all of which have order greater than $48\cdot 23$.) Thus $P$ is
normal in $G$, as required, so our lower bound is attained for $g=24$. We
will now show that this is the least genus for which the bound
is attained.

\medskip

In [{\bf 13}], Kazaz classified the elementary
abelian coverings of the regular hypermaps of genus $2$. In terms of
Fuchsian groups and Riemann surfaces, his results include the following
consequences. Suppose that $\Gamma\geq\Delta\geq K$ where $\Gamma$ is a
triangle group, $\Delta$ is a normal surface subgroup of genus $2$, and
$K$ is a normal subgroup of $\Gamma$ of prime index $p$ in $\Delta$
(so $K$ is a surface group of genus $g=p+1$, and $G=\Gamma/K\leq
{\rm Aut}\,({\cal H}/K)$). Then we have the following possibilities for
$\Gamma,\ Q=\Gamma/\Delta,\ s=|Q|$ and $p$ (all of which occur):
\vskip3pt
\item{a)} $\Gamma=\Gamma(2,8,8),\, Q=C_8,\, s=8,\, p\equiv 1$ mod~$(8)$ or
$p=2$;
\vskip3pt
\item{b)} $\Gamma=\Gamma(4,4,4),\, Q=Q_8,\, s=8,\, p=2$;
\vskip3pt
\item{c)} $\Gamma=\Gamma(2,4,8),\, Q=SD_8,\, s=16,\, p=2$;
\vskip3pt
\item{d)} $\Gamma=\Gamma(5,5,5),\, Q=C_5,\, s=5,\, p\equiv 1$ mod~$(5)$ or
$p=5$;
\vskip3pt
\item{e)} $\Gamma=\Gamma(2,5,10),\, Q=C_{10},\, s=10,\, p\equiv 1$
mod~$(5)$ or $p=5$;
\vskip3pt
\item{f)} $\Gamma=\Gamma(3,6,6),\, Q=C_6,\, s=6,\, p\equiv 1$ mod~$(6)$ or
$p=3$;
\vskip3pt
\item{g)} $\Gamma=\Gamma(2,6,6),\, Q=C_6\times C_2,\, s=12,\, p\equiv 1$
mod~$(6)$ or $p=3$.
\vskip3pt
\noindent(Here $Q_8$ is the quaternion group $\langle a,b \mid a^4=1,\;
a^2=b^2,\;
ba=a^{-1}b\rangle$ of order $8$, and $SD_8$ is the semidihedral group
$\langle a,b \mid a^8=b^2=1,\; ba=a^3b\rangle$ of order $16$.)

These triangle groups $\Gamma$ are all arithmetic, so if $g=p+1$ for any
of the above primes $p$ then $N_{\rm ar}(g)\geq sp>4(g-1)$. Among the genera
$g<24$, those covered by this result are $g=3,\, 4,\, 6,\, 8,\, 12,\, 14,\,
18$ and $20$. To show that each odd $g=2m+1$ satisfies $N_{\rm ar}(g)\geq
6(g-1)$, we can
use a SKE from $\Gamma=\Gamma(2,2,2,6)$ onto $S_3\times D_m$. This leaves
$g=2,\, 10,\, 16$ and $22$ among the genera $g<24$. A SKE $\Gamma(2,3,8)\to
GL(2,3)$ shows that $N_{\rm ar}(2)\geq 48$ (in fact, $N_{\rm
ar}(2)=N(2)=48$). For $g=10$
we can use a SKE from $\Gamma(2,2,2,4)$ onto a split extension of
$C_3\times C_3$ by $D_4$, giving $N_{\rm ar}(10)\geq 72=8(g-1)$. We see in
case (g) that if
$\Gamma=\Gamma(2,6,6)$ then $\Delta$ contains normal surface subgroups $K_3$
and $K_7$ of $\Gamma$, of index $3$ and $7$ in $\Delta$; then
$K=K_3\cap K_7$ is a normal surface subgroup of $\Gamma$ of index $21$ in
$\Delta$, so $N_{\rm ar}(22)\geq 252=12(g-1)$. Finally, for $g=16$ we have a
SKE from $\Gamma(3,3,4)$ onto the alternating group $A_6$, which gives
$N_{\rm ar}(16)\geq 24(16-1)$.

We summarize the results of this section in the following statement:

\medskip

THEOREM 4.1. {\it For all primes $p\equiv 23, 47$ or $59$ {\rm mod}~$(60)$
we have
$$N_{\rm ar}(g)=4(g-1),$$ where $g=p+1$. The least genus $g$
for which $N_{\rm ar}(g)=4(g-1)$ is $g=24$.}

\bigskip\bigskip

\centerline{REFERENCES}

\bigskip

\item{[{\bf 1}]} R.~D.~ACCOLA. On the number of automorphisms of a closed
Riemann surface. {\it Trans.~Amer.~Math.~Soc.} {\bf 131} (1968), 398--408.

\medskip

\item{[{\bf 2}]} R.~D.~ACCOLA. {\it Topics in the Theory of Riemann
Surfaces}. Lecture Notes in Math. {\bf 1595} (Springer-Verlag, 1994).

\medskip

\item{[{\bf 3}]} M.~BELOLIPETSKY. On the number of automorphisms of a
nonarithmetic Riemann surface. {\it Siberian~Math.~J.} {\bf 38} (1997),
860--867.

\medskip

\item{[{\bf 4}]} A.~BOREL. Commensurability classes and volumes of hyperbolic
3-manifolds. {\it Ann. Scuola Norm.~Sup.~Pisa Cl.~Sci.} (4) {\bf 8} (1981),
1--33.

\medskip

\item{[{\bf 5}]} M.~D.~E.~CONDER and R.~KULKARNI. Infinite families of
automorphism groups of Riemann surfaces. In {\it Groups and Geometry}
(eds.~C.~Maclachlan and
W.~J.~Harvey), pp.~47--56. Lond.~Math.~Soc.~Lect.~Note Ser. {\bf 173}
(Cambridge Univ~Press, 1992).

\medskip

\item{[{\bf 6}]} J.~D.~DIXON and B.~MORTIMER. {\it Finite Permutation Groups}.
Grad.~Texts Math. {\bf 163} (Springer-Verlag, 1996).

\medskip

\item{[{\bf 7}]} H.~M.~FARKAS and I.~KRA. {\it Riemann Surfaces}.
Grad.~Texts Math. {\bf 71}
(Springer-Verlag, 1980).

\medskip

\item{[{\bf 8}]} B.~HUPPERT. {\it Endliche Gruppen I\/} (Springer-Verlag, 1967).

\medskip

\item{[{\bf 9}]} B.~HUPPERT and N.~BLACKBURN. {\it Finite Groups III\/}
(Springer-Verlag, 1982).

\medskip

\item{[{\bf 10}]} A.~HURWITZ. \"Uber algebraische Gebilde mit eindeutige
Transformationen in sich. {\it Math.~Ann.} {\bf 41} (1893), 403--442.

\medskip

\item{[{\bf 11}]} G.~A.~JONES and D.~SINGERMAN. {\it Complex Function Theory: an
Algebraic and Geometric Viewpoint} (Cambridge University Press, 1987).

\medskip

\item{[{\bf 12}]} S.~KATOK. {\it Fuchsian Groups} (University of Chicago Press,
1992).

\medskip

\item{[{\bf 13}]} M.~KAZAZ. Finite Groups and Surface Coverings.
Ph.D.~thesis. University of
Southampton (1997).

\medskip

\item{[{\bf 14}]} C.~MACLACHLAN. A bound for the number of automorphisms
of a compact Riemann surface. {\it J.~London Math.~Soc.} {\bf 44}
(1968), 265--272.

\medskip

\item{[{\bf 15}]} C.~MACLACHLAN and G.~ROSENBERGER. Two-generator arithmetic
Fuchsian groups, II. {\it Math.~Proc.~Camb.~Phil.~Soc.} {\bf 111} (1992),
7--24.

\medskip

\item{[{\bf 16}]} C.-H.~SAH. Groups related to compact Riemann surfaces.
{\it Acta Math.} {\bf 123} (1969), 13--42.

\medskip

\item{[{\bf 17}]} H.~A.~SCHWARZ. Ueber diejenigen algebraischen Gleichungen
zwischen zwei ver-\"anderlichen Gr\"ossen, welche eine Schaar rationaler,
eindeutig umkehrbarer Transformationen in sich selbst zulassen.
{\it J.~reine angew.~Math.} {\bf 87}, 139--145; {\it Ges.~Math. Abh.~II\/},
pp.~285--291
(1890, reprinted Chelsea, 1972).

\medskip

\item{[{\bf 18}]} K.~TAKEUCHI. Arithmetic triangle groups.
{\it J.~ Math.~Soc.~Japan} {\bf 29} (1977), 91--106.

\medskip

\item{[{\bf 19}]} K.~TAKEUCHI. Arithmetic Fuchsian groups with signature
$(1;e)$. {\it J.~Math.~Soc.~Japan} {\bf 35} (1983), 381--407.

\medskip

\item{[{\bf 20}]} M.-F.~VIGN\'{E}RAS. {\it Arithm\'{e}tique des Alg\`{e}bres de
Quaternions}. Lecture Notes in Math. {\bf 800} (Springer-Verlag, 1980).

\vfil\eject

\centerline{Table: Signatures $\sigma\in\Sigma$}

\bigskip

\centerline{
\vbox{\tabskip=0pt \offinterlineskip
\def\tablerule{\noalign{\hrule}}
\halign to180pt{\strut#& \vrule#\tabskip=1em plus2em&
  #& \vrule#& \hfil#\hfil& \vrule#&
  #\hfil& \vrule#\tabskip=0pt\cr\tablerule
&&\omit\hidewidth $\sigma$\hidewidth&&
\omit\hidewidth $\mu(\Gamma)$\hidewidth&&
\omit\hidewidth $s/r$\hidewidth&\cr\tablerule
&&(2,3,7)&&$\pi/21$&&84&\cr\tablerule
&&(2,3,8)&&$\pi/12$&&48&\cr\tablerule
&&(2,3,9)&&$\pi/9$&&36&\cr\tablerule
&&(2,3,10)&&$2\pi/15$&&30&\cr\tablerule
&&(2,3,11)&&$5\pi/33$&&132/5&\cr\tablerule
&&(2,3,12)&&$\pi/6$&&24&\cr\tablerule
&&(2,3,14)&&$4\pi/21$&&21&\cr\tablerule
&&(2,3,16)&&$5\pi/24$&&96/5&\cr\tablerule
&&(2,3,18)&&$2\pi/9$&&18&\cr\tablerule
&&(2,3,24)&&$\pi/4$&&16&\cr\tablerule
&&(2,3,30)&&$4\pi/15$&&15&\cr\tablerule
&&(2,4,5)&&$\pi/10$&&40&\cr\tablerule
&&(2,4,6)&&$\pi/6$&&24&\cr\tablerule
&&(2,4,7)&&$3\pi/14$&&56/3&\cr\tablerule
&&(2,4,8)&&$\pi/4$&&16&\cr\tablerule
&&(2,4,10)&&$3\pi/10$&&40/3&\cr\tablerule
&&(2,4,12)&&$\pi/3$&&12&\cr\tablerule
&&(2,4,18)&&$7\pi/18$&&72/7&\cr\tablerule
&&(2,5,5)&&$\pi/5$&&20&\cr\tablerule
&&(2,5,6)&&$4\pi/15$&&15&\cr\tablerule
&&(2,5,8)&&$7\pi/20$&&80/7&\cr\tablerule
&&(2,5,10)&&$2\pi/5$&&10&\cr\tablerule
&&(2,5,20)&&$\pi/2$&&8&\cr\tablerule
&&(2,5,30)&&$8\pi/15$&&15/2&\cr\tablerule
&&(2,6,6)&&$\pi/3$&&12&\cr\tablerule
&&(2,6,8)&&$5\pi/12$&&48/5&\cr\tablerule
&&(2,6,12)&&$\pi/2$&&8&\cr\tablerule
&&(2,7,7)&&$3\pi/7$&&28/3&\cr\tablerule
&&(2,7,14)&&$4\pi/7$&&7&\cr\tablerule
&&(2,8,8)&&$\pi/2$&&8&\cr\tablerule
&&(2,8,16)&&$5\pi/8$&&32/5&\cr\tablerule
&&(2,9,18)&&$2\pi/3$&&6&\cr\tablerule
&&(2,10,10)&&$3\pi/5$&&20/3&\cr\tablerule
&&(2,12,12)&&$2\pi/3$&&6&\cr\tablerule
&&(2,12,24)&&$3\pi/4$&&16/3&\cr\tablerule
&&(2,15,30)&&$4\pi/5$&&5&\cr\tablerule
&&(2,18,18)&&$7\pi/9$&&36/7&\cr\tablerule
\noalign{\smallskip}\cr}}
\quad
\vbox{\tabskip=0pt \offinterlineskip
\def\tablerule{\noalign{\hrule}}
\halign to180pt{\strut#& \vrule#\tabskip=1em plus2em&
  #& \vrule#& \hfil#\hfil& \vrule#&
  #\hfil& \vrule#\tabskip=0pt\cr\tablerule
&&\omit\hidewidth $\sigma$\hidewidth&&
\omit\hidewidth $\mu(\Gamma)$\hidewidth&&
\omit\hidewidth $s/r$\hidewidth&\cr\tablerule
&&(3,3,4)&&$\pi/6$&&24&\cr\tablerule
&&(3,3,5)&&$4\pi/15$&&15&\cr\tablerule
&&(3,3,6)&&$\pi/3$&&12&\cr\tablerule
&&(3,3,7)&&$8\pi/21$&&21/2&\cr\tablerule
&&(3,3,8)&&$5\pi/12$&&48/5&\cr\tablerule
&&(3,3,9)&&$4\pi/9$&&9&\cr\tablerule
&&(3,3,12)&&$\pi/2$&&8&\cr\tablerule
&&(3,3,15)&&$8\pi/15$&&15/2&\cr\tablerule
&&(3,4,4)&&$\pi/3$&&12&\cr\tablerule
&&(3,4,6)&&$\pi/2$&&8&\cr\tablerule
&&(3,4,12)&&$2\pi/3$&&6&\cr\tablerule
&&(3,5,5)&&$8\pi/15$&&15/2&\cr\tablerule
&&(3,6,6)&&$2\pi/3$&&6&\cr\tablerule
&&(3,6,18)&&$8\pi/9$&&9/2&\cr\tablerule
&&(3,8,8)&&$5\pi/6$&&24/5&\cr\tablerule
&&(4,4,4)&&$\pi/2$&&8&\cr\tablerule
&&(4,4,5)&&$3\pi/5$&&20/3&\cr\tablerule
&&(4,4,6)&&$2\pi/3$&&6&\cr\tablerule
&&(4,4,9)&&$7\pi/9$&&36/7&\cr\tablerule
&&(4,5,5)&&$7\pi/10$&&40/7&\cr\tablerule
&&(4,6,6)&&$5\pi/6$&&24/5&\cr\tablerule
&&(5,5,5)&&$4\pi/5$&&5&\cr\tablerule
&&(2,2,2,3)&&$\pi/3$&&12&\cr\tablerule
&&(2,2,2,4)&&$\pi/2$&&8&\cr\tablerule
&&(2,2,2,5)&&$3\pi/5$&&20/3&\cr\tablerule
&&(2,2,2,6)&&$2\pi/3$&&6&\cr\tablerule
&&(2,2,2,7)&&$5\pi/7$&&28/5&\cr\tablerule
&&(2,2,2,8)&&$3\pi/4$&&16/3&\cr\tablerule
&&(2,2,2,9)&&$7\pi/9$&&36/7&\cr\tablerule
&&(2,2,2,10)&&$4\pi/5$&&5&\cr\tablerule
&&(2,2,2,12)&&$5\pi/6$&&24/5&\cr\tablerule
&&(2,2,2,14)&&$6\pi/7$&&14/3&\cr\tablerule
&&(2,2,2,18)&&$8\pi/9$&&9/2&\cr\tablerule
&&(2,2,2,22)&&$10\pi/11$&&22/5&\cr\tablerule
&&(2,2,3,3)&&$2\pi/3$&&6&\cr\tablerule
&&(2,2,3,4)&&$5\pi/6$&&24/5&\cr\tablerule
&&(2,2,3,5)&&$14\pi/15$&&30/7&\cr\tablerule
\noalign{\smallskip}\cr}}
}

\end